\documentclass[a4paper,twoside,12pt]{amsart}
\usepackage{amsfonts}
\usepackage{mathrsfs}
\usepackage{amsmath,graphics}
\usepackage{amssymb}
\usepackage{indentfirst,latexsym,bm,amsmath,pstricks,amssymb,amsthm,graphicx,fancyhdr,float,fullpage}
\usepackage[all]{xy}
\usepackage{amsthm}
\usepackage{amscd}
\usepackage[CJKbookmarks=true]{hyperref}
\usepackage{color}
\usepackage{chngpage}
\usepackage{aliascnt}
\usepackage{hyperref}
\usepackage{appendix}
\usepackage{soul}

\usepackage{url}

\newtheorem{theorem}{Theorem}[section]

\newtheorem*{acknowledgement*}{\protect \acknowledgementname}

\newaliascnt{setup}{theorem}

\aliascntresetthe{setup}

\newaliascnt{question}{theorem}

\aliascntresetthe{question}

\newaliascnt{lemma}{theorem}

\aliascntresetthe{lemma}

\newaliascnt{assumption}{theorem}

\aliascntresetthe{assumption}

\newaliascnt{conjecture}{theorem}
\newtheorem{conjecture}[conjecture]{Conjecture}
\aliascntresetthe{conjecture}

\newaliascnt{proposition}{theorem}

\aliascntresetthe{proposition}

\newaliascnt{corollary}{theorem}
\newtheorem{corollary}[corollary]{Corollary}
\aliascntresetthe{corollary}

\newaliascnt{problem}{theorem}

\aliascntresetthe{problem}

\newaliascnt{claim}{theorem}

\aliascntresetthe{claim}

\theoremstyle{definition}

\newaliascnt{definition}{theorem}

\aliascntresetthe{definition}

\newaliascnt{example}{theorem}

\aliascntresetthe{example}

\theoremstyle{remark}

\newaliascnt{remark}{theorem}

\aliascntresetthe{remark}

\newaliascnt{remarks}{theorem}

\aliascntresetthe{remarks}

\def \({$($}
\def \){$)$}

\providecommand{\acknowledgementname}{Acknowledgement}

\def \GL2{\mathrm{GL}_2}
\def \SL2{\mathrm{SL}_2}

\numberwithin{equation}{section}

\begin{document}
\title{Constructing algebraic solutions of Painleve VI equation from $p$-adic Hodge theory and Langlands Correspondence}
\author{Jinbang Yang}
\email{yjb@mail.ustc.edu.cn}
\address{School of Mathematical Sciences, University of Science and Technology of China, Hefei, Anhui 230026, PR China}
\author{Kang Zuo}
\email{zuok@uni-mainz.de}
\address{Institut f \"ur Mathematik, Universit \"at Mainz, Mainz 55099, Germany}

\maketitle
\begin{abstract} We construct infinitely many non-isotrivial families of abelian varieties over given four punctured projective lines. These families lead to algebraic solutions of Painleve VI equation. Finally, based on a recent paper by Lin-Sheng-Wang, we prove  a complete characterization for the locus of motivic Higgs bundles in the moduli space as fixed points of an ``additive'' self-map.
This is a note based on the lecture given by the second named author on 04.Nov.2022 at Tsinghua University.	
\end{abstract}

\section{\bf Introduction}
Let $R$ be a commutative ring with identity and $X$ be a scheme over $R$. An abelian scheme $A$ over $X$ together with a polarization $\mu$ is called of \emph{of $\mathrm{GL}_2$-type}, if there exists a number field $K$ of degree $\dim_X A$ such that the ring of integers $\mathcal O_K$ can be embedded into the endomorphism ring $\mathrm{End}_{\mu}(A/X)$. We will call the abelian scheme \emph{of $\mathrm{GL}_2(K)$-type}, if we want to emphasize the role of $K$. 

Let $f\colon A\rightarrow X$ be an abelian scheme of $\mathrm{GL}_2(K)$-type. Let $D$ denote the discriminant locus and let $X^0$ denote the complement of $D$ in $X$. Let $\Delta\subset A$ denote the inverse image of $D$ under the structure morphism $f$ and let $A^0$ denote the complement of $\Delta$ in $A$. Then we obtain the smooth abelian scheme $f^0\colon A^0\rightarrow X^0$.
\begin{equation}
	\xymatrix{
		A^0 \ar[r] \ar[d]^{f^0} & A \ar[d]^f & \Delta \ar[d]\ar[l] \\
		X^0 \ar[r] & X & D \ar[l] \\
	}
\end{equation} 

For $R= \mathbb C$ we consider
the Betti-local system
\[\mathbb V= R^1_\mathrm{B}f^0_* \mathbb Z_{A^0}\] 
attached to
$f^0$, which is a $\mathbb{Z}$-local system over the base $X^0$. Since $f$ is of $\mathrm{GL}_2(K)$-type, the action of $\mathcal O_K$ on $f$ induces an action of $K$ on the $\overline{\mathbb Q}$-local system $\mathbb V \otimes \overline{\mathbb Q}$. Taking the $K$-eigen sheaves decomposition 
\[\mathbb V \otimes \overline{\mathbb Q}= \bigoplus_{i=1}^g \mathbb L_i.\]
Then there $\mathbb L_i$'s are of rank-$2$ over $X^0$ and defined over the ring of integers of some number field. On the other hand, consider the logarithmic de Rham bundle attached to the abelian scheme $f$ and denote
\[(V,\nabla)= R^1_\mathrm{dR}f_*\Big(\Omega^*_{A/X}(\log \Delta),\mathrm{d}\Big).\]
On this de Rham bundle, there is a canonical filtration satisfying Griffiths transversality given by 
\[E^{1,0}:=R^0f_* \Omega^1_{A/X}(\log \Delta))\subset V.\]
Taking the grading with respect to this filtration, one gets a logarithmic graded Higgs bundle, which is so-called Kodaira-Spencer map attached to $f$
\begin{equation} \label{equ:KS}
	(E,\theta):=(E^{1,0} \oplus E^{0,1},\theta):= \mathrm{Gr}_{E^{1,0}} (V,\nabla) = \Big(R^0f_*\Omega^1_{A/X}(\log \Delta)\oplus R^1f_*\mathcal O_A,\mathrm{Gr}(\nabla)\Big).
\end{equation}
Since $f$ is of $\mathrm{GL}_2(K)$-type, one also gets a $K$-eigen decomposition of the Higgs bundle
\begin{equation} \label{equ:decomp_Higgs}
	(E,\theta)=\bigoplus_{i=1}^{g} (E,\theta)_i.
\end{equation}
Under Hitchin-Simpson's non-abelian Hodge theory, these eigensheaves $\{(E,\theta)_i\}_{i=1,\cdots,g}$ are just those Higgs bundles correspond to the local systems $\{\mathbb L_i\}_{i=1,\cdots,g}$.

Local systems and Higgs bundles arising as subquotient of local systems and Higgs bundles of  some families of varieties are called \emph{motivic}. Sometimes they are also called \emph{coming from geometry origin}. Simpson had found a characterization for a rank-2 local system being motivic:
\begin{theorem}[Simpson]
	A rank-2 local system $\mathbb L$ is motivic if and only if the following two conditions hold:
	\begin{enumerate}
		\item $\mathbb L$ is defined over the ring of integers of some number field, and
		\item for each element $\sigma \in \mathrm{Gal}(\overline{\mathbb Q}/\mathbb Q)$
		the Higgs bundle corresponding to the Galois conjugation $\mathbb L^\sigma$ is again graded.
	\end{enumerate}
\end{theorem}

\begin{conjecture}[Simpson]
A rigid local system is motivic.
\end{conjecture}
Simpson's conjecture has been proved for the case of rank-2 by Corlette-Simpson \cite{CS} and rank-3 by Langer-Simpson \cite{LS} for cohomologically rigid local systems. The conjecture predicts that any rigid local system $\mathbb L$ shall enjoy all properties of motivic local systems. For example,
\begin{itemize}
		\item  its corresponding filtered de Rham bundle is isomorphic to the underlying filtered de Rham bundle of some Fontaine-Faltings modules at almost all places, and  
		\item  if $\mathbb L$   is in addtion cohomologically rigid, then it is defined over the ring of integers of some number field.

\end{itemize}
Those two properties have been verified by Esnault-Groechenig recently \cite{EG20}, \cite{EG18}.

In this note, we take $X$ as the complex projective line $\mathbb P^1$ and $D$ as the $4$ punctures $\{0,1,\infty,\lambda\}$. Our goal is find some motivic rank-$2$ logarithmic graded Higgs bundles over
$(\mathbb P^1,\{0,1,\infty,\lambda\})$, which are coming from $\mathrm{GL}_2$-type abelian schemes over $\mathbb P^1$ with discriminant locus contained in $\{0,1,\infty,\lambda\}$.

Beavuille has shown that there exist exactly 6 non-isotrivial families of elliptic curves over $\mathbb P^1$ with semistable reductions over $\{0,1,\infty,\lambda_i\}$ for ${1\leq i\leq 6}$. All of them are modular curves with respect to certain mixed level structures. Based on Beavuille's result, Viehweg-Zuo have shown that there no more non-isotrivial abelian schemes over $\mathbb P^1$ of $\mathrm{GL}_2(K)$-type with $[K:\mathbb Q]\geq2$ and with semistable reductions over $\{0,1,\infty,\lambda\}$. So except Beavuille's example any non-isotrivial smooth abelian schemes over $\mathbb P^1 \setminus \{0,1,\infty,\lambda\}$ of $\mathrm{GL}_2(E)$-type must have non-semistable reduction at some point in $\{0,1,\infty,\lambda\}$. In this case, the some eigenvalues of the local monodromies of motivic local system must be roots of unity other than $1$. 

In this note, we consider the simplest situation: motivic rank-2 local systems whose local monodromies around $\{0,1,\infty,\lambda\}$ are unipotent and around $\infty$ is quasiunipotent with eigenvalues $\{-1,-1\}$. We call such type local monodromy to be \emph{of type-$(1/2)_\infty$}, or just \emph{of type-$(1/2)$}.

\begin{theorem}\label{thm:main}
For given $\lambda\in \mathbb P^1\setminus\{0,1,\infty\}$, there exists infinitely many non-isotrivial abelian schemes of $GL_2$-type over $\mathbb P^1\setminus\{0,1,\infty,\lambda\}$ with the associated rank-$2$ eigen local systems being of type-$(1/2)_\infty$.
\end{theorem}

Let $M_{0,n}$ denote the moduli space of $n$-punctured projective lines and let $S_{0,n}$ denote the total space of the universal family of $n$-punctured projective lines with structure morphism
\[p_n\colon S_{0,n}\rightarrow M_{0,n}.\]
Then $M_{0,4}\simeq \mathbb P^1\setminus \{0,1,\infty\}$, and $S_{0,4}=\bigcup\limits_{\lambda \in M_{0,4}}\big(\mathbb P^1\setminus\{0,1,\infty,\lambda\}\big)$ is an algebraic surface. Once we vary the parameter $\lambda$, \autoref{thm:main} implies the following result: 
\begin{theorem} \label{thm:mainII}
There exist infinitely many  non-isotrivial abelian schemes 
	$$f: A\to \widetilde S_{0,4}$$
of $\mathrm{GL}_2$-type over fiber products of finite \'etale base changes
\begin{equation*}
\xymatrix@C=2cm{ 
 \widetilde{S}_{0,4}:=\widetilde {M}_{0,4} \times_{M_{0,4}} S_{0,4} \ar[r] \ar[d] & {S}_{0,4} \ar[d]^{p_4}\\
 \widetilde {M}_{0,4}^{} \ar[r]^{\text{finite \'etale}} & M_{0,4}\\	
} 
\end{equation*} 
and such that the local monodromies of $f$ around $\{0,1,\lambda\}$ are unipotent and around $\{\infty\}$ is quasi-unipotent with all eigenvalues being $-1$.
\end{theorem}

\begin{corollary}[Corollary 4.] Let $f: A\to \widetilde S_{0,4}$ be a family given in \autoref{thm:mainII}. Then all rank-$2$ eigen local systems associated to the family $f$ are algebraic solutions of Painleve VI equation of the type-$(1/2)_\infty$.
\end{corollary}

For given $\lambda\in \mathbb P^1\setminus\{0,1,\infty\}$, any family $f_\lambda\colon A_\lambda\rightarrow \mathbb P^1$ in \autoref{thm:main} has semistable reduction over $\{0,1,\lambda\}$ and potentially semistable reduction over $\infty$. Thus the eigen Higgs bundles $(E,\theta)_i$ associated to this family (constructed in \ref{equ:decomp_Higgs}) have the following form
\begin{equation} \label{equ:Higgs_form}
	E_i=\mathcal O\oplus \mathcal O(-1), \qquad \theta_i\colon \mathcal O\xrightarrow{\neq0} \mathcal O(-1) \otimes \Omega^1_{\mathbb P^1}(\log \{0,1,\infty,\lambda \})
\end{equation}
and are endowed with natural parabolic structures on the punctures $\{0,1,\infty,\lambda\}$ of type-$(1/2)_\infty$.
Here type-$(1/2)_\infty$ parabolic structures means that the parabolic structures at $0$, $1$ and $\lambda$ are trivial and the parabolic filtration at $\infty$ is 
\[\left(E_{i}\mid_\infty\right)_\alpha=\left\{\begin{array}{cc}
E_{i}\mid_\infty & 0\leq\alpha\leq1/2,\\
0 & 1/2 <\alpha < 1.\\
\end{array}\right.\]
Let $M^{gr {1\over 2}}_{Hig\lambda} $ denote the moduli space of rank-2 semi-stable graded Higgs bundles over $\mathbb P^1$ with the parabolic structure on $\{0,1,\infty,\lambda \}$ of type-$(1/2)_\infty$. Then any Higgs bundle $(E,\theta)\in M^{gr {1\over 2}}_{Hig\lambda} $ is parabolic stable and has the form as in \eqref{equ:Higgs_form}. 

In view of $p$-adic Hodge theory, a Higgs bundle $(E,\theta) $ over the Witt ring $W(\mathbb F_q)$ realized by an abelian scheme over $W(\mathbb F_q)$ of $\text{GL}(E)_2$-type has to be the grading of an $K$-eigen sheaf of the Fontaine-Faltings module 
attached to the abelian scheme. Hence, by Lan-Sheng-Zuo functor $(E,\theta)$ is \emph{periodic} on $M^{gr {1\over 2}}_{Hig\lambda} $ over $W(\mathbb F_q)$.

One identifies 
 $$M^{gr {1\over 2}}_{Hig\lambda} = \mathbb P^1$$ 
 by sending $(E,\theta)$ to the zero locus of the Higgs map $(\theta)_0\in \mathbb P^1$, and takes then the elliptic curve $C_{\lambda}$ of the Weierstrass form $y^2=z(z-1)(z- \lambda)$ as the double cover
 $$\pi: C_\lambda\to \mathbb P^1$$ ramified on $\{0,1,\infty,\lambda \}$. 
 \begin{conjecture} [Sun-Yang-Zuo \cite{SYZ}] \label{conj:SYZ}
	 The self-map induced by Higgs-de Rham flow on $ {M^{gr {1\over 2}}_{Hig\lambda}}$ over ${\mathbb F_q} $ comes from multiplication map by $p$ on the associated elliptic curve over $\mathbb F_q$
\[\xymatrix{
 		& C_\lambda \ar[d] \ar[r]^{[p]} & C_\lambda \ar[d] & \\
 		M_{Higg,\lambda}^{gr}\ar@/_12pt/[rrr]_{\phi} \ar@{=}[r] & \mathbb P^1 \ar[r] & \mathbb P^1 \ar@{=}[r] & M_{Higg,\lambda}^{gr} \\
 	}\]
\end{conjecture}

The conjecture implies two things:
\begin{enumerate}
	\item a Higgs bundle $(E,\theta)$ is periodic if and only if $\pi^{-1}(\theta)_0$ is a torsion point in $C_\lambda$ and of order $p^f-1.$
	\item for a prime $p>2$ and assume $C_\lambda$ is supersingular then $\phi_\lambda(z)=z^{p^2}$. Hence, any Higgs bundle $(E,\theta)\in { M^{gr {1\over 2}}_{Hig\lambda}}(\overline{\mathbb F}_{q})$ is periodic.
\end{enumerate} 

The Conjecture has been checked by Sun-Yang-Zuo for $p<50.$ Very recently it has been proved by 
 Lin-Sheng-Wang and becomes a theorem.
 \begin{theorem} [Lin-Sheng-Wang \cite{LSW}] \label{thm:LSW}
 	\autoref{conj:SYZ} holds true.
 \end{theorem}

The technique used in \autoref{thm:mainII} combined with \autoref{thm:LSW} lead us to prove the following result on a characterization for motivic Higgs bundles contained in $M^{gr {1\over 2}}_{Hig\lambda}$.
 \begin{theorem} \label{thm:torsion}
A Higgs bundle $(E,\theta)\in M^{gr {1\over 2}}_{Hig\lambda} $ 
is an eigensheaf of the Kodaira-Spencer map attached to an abelian scheme of $\text{GL}_2$-type if and only if
$\pi^{-1}(\theta)_0$ is a torsion point in $C_\lambda$.
\end{theorem}

\section{ Discussion on \autoref{thm:main}}

The  underlying principle behind \autoref{thm:main} is very simple, the so-called isomonodromy deformation for motivic local systems over mixed characteristic.
To illustrate the idea,  let's first look at the situation over complex numbers. We assume, there exists an abelian scheme 
$f_{\lambda_0}: A_{\lambda_0}\to \mathbb P^1$ of $\text{GL}_2(K)$-type over complex numbers with bad reduction on $\{0,\,1,\,\lambda_0,\,\infty\}$ of type-(1/2).
Then the filtered logarithmic de Rham bundle decomposes as $K$-eigen sheaves
$$(V, \nabla, E^{1,0})=:( R^1_ \mathrm{dr}f_* \Omega^*_{A_{\lambda_0}/\mathbb P^1}(\log \Delta), d), \, R^0f_* \Omega^1_{A_{\lambda_0}/\mathbb P^1}(\log \Delta))=\bigoplus_{i=1}^g(V,\nabla, E^{1,0})_i,$$
where each eigen sheaf has the form
$$ (V,\nabla, E^{1,0})_i\simeq (\mathcal O\oplus \mathcal O(-1), \nabla_i, \mathcal O). $$ 
\\
Consider a family of 4-punctured projective line
$$(\mathbb P^1, \{0,\,1,\,\lambda,\, \infty\})_{\hat U_{\lambda_0}}\to \hat U_{\lambda_0}$$
over a formal neighborhood $\hat U_{\lambda_0}\subset M_{0\, 4}$ of $\lambda_0.$
Then by forgetting the Hodge filtration the de Rham bundle extends to a de Rham bundle 
$(V, \nabla)_{\hat U_{\lambda_0}}$
over $(\mathbb P^1, \{0,\,1,\,\lambda,\, \infty\})_{\hat U_{\lambda_0}}$. It is known the abelian scheme extends over  $(\mathbb P^1, \{0,\,1,\,\lambda,\, \infty\})_{\hat U_{\lambda_0}}$
if and only if the Hodge filtration $E^{1,0}$ extends to a sub bundle in the de Rham bundle $(V, \nabla)_{\hat U_{\lambda_0}}$. Using the $K$-eigen sheave decomposition we see that the obstruction for extending the Hodge filtration $E^{1,0}=\bigoplus_{i=1}^g \mathcal O$
lies in 
$\bigoplus_{i=1}^g H^1(\mathbb P^1, \mathcal O(-1))=0$. Hence, the abelian scheme $f_{\lambda_0}$
extends over the base $(\mathbb P^1, \{0,\,1,\,\lambda,\, \infty\})_{\hat U_{\lambda_0}}$.

Back to the situation  over mixed characteristic, along the diagram below. One 
\begin{itemize}
	\item starts with moduli space ${ M^{gr {1\over 2}}_{Hig\lambda}}$ of rank-$2$ stable graded Higgs bundles on $\mathbb P^1_{\mathbb Z_q}$ of parabolic type-$(1/2)_\infty$ on $\{0,1,\lambda,\infty\}$,
	\item finds periodic Higgs bundles over $\mathbb P^1_{\mathbb F_q}$(i.e. fixed points of iterations of the self map on ${ M^{gr {1\over 2}}_{Hig\lambda}}\otimes \mathbb F_q$ induced by Higgs-de Rham flow),
	\item gets Fontaine-Faltings modules via Lan-Sheng-Zuo functor from those periodic Higgs bundles and lifts these modules to $\mathbb P^1_{\mathbb Z_q}$, 
	\item obtains rank-2 $\ell$-adic local systems on $\mathbb P^1-\{0,\, 1,\,\lambda,\,\infty\}$ over $\mathbb F_q$ by forgetting Hodge filtration in the Fontaine-Faltings module, tensoring with $\mathbb Q_p$, and applying Deligne's $p$-$\ell$ companion conjecture proven by Abe \cite{Ab}, 
	\item finds $\text{GL}_2$-type abelian schemes over $\mathbb P^1_{\mathbb F_q}$ with bad reductions of type-$(1/2)_\infty$ realizing those $\ell$-adic local systems via Drinfeld's theorem on Langlands correspondence. 
	\item modifies those $\mathrm{GL}_2$-type abelian schemes up to some $p$-isogeny and lifts these modified $\mathrm{GL}_2$-type abelian schemes to $\mathbb P^1_{\mathbb Z_q}$. 
\end{itemize}  
\[\xymatrix@C=3cm @ R=2cm{
	M_{\rm p-adic}^{\rm cris}/\chi
	\ar@{<->}[r]^-{\text{Fontaine-Faltings}}_-{p\text{\rm-adic RH}} 
	& M_{\rm dR}^{FF}/\chi
	\ar@{^(->}[r]^{\text{forgetting Hodge}}_{\text{filtration, }\otimes \mathbb Q_p} 
	\ar@{<->}[d]_{\text{Higgs-de Rham flow}}^{\text{by Lan-Sheng-Zuo}} 
	& M^{\rm F-isoc}/\chi
	\ar@{<->}[r]^{\text{Deligne's } p-\ell}_{\text{companion by Abe}} 
	& M_{\ell-\text{adic}}^{}\\
	& M_{\rm Higg}^{\rm per} \ar@{^(->}[r]
	& { M^{\rm gr {1\over 2}}_{\rm Hig\,\lambda}} 
	& \\
}\]

And a type of  boundedness and rigidity arguments of classifying maps from the log base curve $(\mathbb P^1,\{0,1,\lambda,\infty\})$
into the fine moduli space of polarized abelian varieties shows that the abelian scheme lifts over complex numbers for any $\lambda \in M_{0,4} (\mathbb C).$

Below we give more detailed explanations on the technique issue in  each step:
\subsection{\bf Constructing Fontaine-Faltings modules over $\mathbb Z_q$ from semistable Higgs bundles over $\mathbb F_q$ via Higgs de Rham flow.}

Take a prime $p$ and an element $\lambda \in \mathbb Z_{q^2}$
such that the modulo $p$ reduction of $C_\lambda$ is supersingular. Consider the moduli space $M^{gr {1\over 2}}_{Hig\lambda}\otimes \mathbb F_{q^2} $ of rank-2 stable logarithmic graded Higgs bundles over $(\mathbb P^1, D)_{\mathbb F_{q^2}}$ with parabolic structures on $ D$ of type-$(1/2)$. Identifying $M^{gr {1\over 2}}_{Hig\lambda}\otimes \mathbb F_{q^2} $ with the projective line $\mathbb P^1_{\mathbb F_{q^2}}$, the self-map induced by Higgs-de Rham flow is given by 
$z\mapsto z^{p^2}$. We show that, for any $n\geq1$, each $(\overline E,\overline \theta)\in M^{gr {1\over 2}}_{Hig\lambda}(\overline{\mathbb F}_{q^{2n}}) $ is automatically periodic and lifts uniquely to
a periodic Higgs bundle $(E,\theta)\in M^{gr {1\over 2}}_{Hig\lambda}( \mathbb Z_{q^{2n}}).$ Hence, the Lan-Sheng-Zuo functor
induces following bijection
$$ M^{gr{1\over 2}}_{Hig\lambda} (\mathbb F_{q^{2n}} ) \simeq M^{FF\,{1\over 2}}_{dR\,\lambda}(\mathbb Z_{q^{2n}})/ \chi,$$
where the right hand side is $M^{FF\,{1\over 2}}_{dR\,\lambda}(\mathbb Z_{q^{2n}})$ modulo an equivalent relation. The $M^{FF\,{1\over 2}}_{dR\,\lambda}(\mathbb Z_{q^{2n}})$ is the set of rank-$2$ Fontaine-Faltings module over $(\mathbb P^1,D)_{\mathbb Z_{q^{2n}}}$ such that the monodromy is of type-$(1/2)$, and two Fontaine-Faltings modules in $M^{FF\,{1\over 2}}_{dR\,\lambda}(\mathbb Z_{q^{2n}})$ are called equivalent if they are differed by a constant rank-$1$ Fontaine-Faltings modules. 

\subsection{\bf Deligne's $p$-$\ell$ companion conjecture solved by Abe.}

Consider the functors defined by forgetting Hodge filtration and tensoring with $\mathbb Q_p$, which induce an injective map between a set of logarithmic Fontaine-Faltings modules
with parabolic structure-(${1\over 2}$) and a set of logarithmic $F$-isocrystals. These functors keeps the same type parabolic structure, and the result logarithmic $F$-isocrystals have a fixed constant determinant, the $F$-isocrystal $\mathcal{E}_{\rm cy}$ associated to the cyclotomic character. The injectivity follows the fact that the underlying log-parabolic de Rham bundles are parabolic stable mod $p$.

We apply furtherly the forgetful functor by sending logarithmic $F$-isocrystals to overconvergent F-isocrystals, which is known injective by work of Kedlaya \cite[Proposition 6.3.2]{Ked07}. Hence, by putting every thing together we obtain an injective map 
$$ M^{gr{1\over 2}}_{Hig\lambda} (\mathbb F_{q^{2n}} ) \simeq M^{FF\,{1\over 2}}_{dR\,\lambda}(\mathbb Z_{q^{2n}})/ \chi\hookrightarrow M^{F-iso\,{1\over 2}\, \dagger}_\lambda(\mathbb Q_{q^{2n}})/\chi ,$$
where $ M^{F-iso\,{1\over 2}\dagger}_\lambda(\mathbb Q_{q^{2n}})$ is the set of overconvergent F-isocrystals with parabolic structure of type-$({1/2})$ and determinant
$\mathcal{E}_{\rm cy}$.

We choose a prime $\ell\not=p$ and fix an isomorphism $\phi: \overline {\mathbb Q}_p\simeq \overline {\mathbb Q}_\ell$. 
Given an overconvergent $F$-isocrystal $(V,\nabla,\Phi)^\dagger\in M^{F-iso\,{1\over 2}\dagger}_\lambda(\mathbb Q_{q^{2n}})$, by applying 
Deligne's $p-\ell$ companion proved by Abe to $(V,\nabla,\Phi)^\dagger$
we find a rank-2 $\ell$-adic irreducible local system $\mathbb L$ with cyclotomic determinant corresponding to $(V,\nabla,\Phi)^\dagger$ in the sense the 
characteristic polynomial of $\Phi_x$ on $V_x$ and the characteristic polynomial of $\sigma_x$ on $\mathbb L_x$ via $\phi$ at all close point $x$ in the close fiber. By the compatibility of local-global Langlands correspondence we 
find the local monodromies of $\mathbb L$ around $\{0,1,\lambda\}$ are unipotent and around $\infty$ is quasi-unipotent with eigenvalues $\{-1,\,-1\}.$ 
Hence, we obtain a bijective functor
 $$ M^{F-iso\, {1\over 2}\dagger}_\lambda(\mathbb Q_{q^{2n}})/\chi \simeq 
 M^{{1\over 2}}_{\ell-\text{adic}\,\lambda}(\mathbb F_{q^{2n}})/\sim,$$
where $ M^{ {1\over 2}}_{\ell-\text{adic}\,\lambda}(\mathbb F_{q^{2n}})$ is the set of equivalent classes of rank-2 $\ell$-adic local systems
over $(\mathbb P^1- D)_{\mathbb F_q^{2n}}$ with cyclotomic determinant, unipotent local monodromies around $\{0,1,\lambda\}$ and quasi-unipotent local monodromy around $\infty$ of eigenvalues $\{-1,-1\}$. The $\sim$ stands for an equivalent relation, two local systems $\mathbb L$ and $\mathbb L'$ are called equivalent if the restriction of them on geometric fundamental group are isomorphic. As we show that the restriction of $\mathbb L$ to the geometric fundamental group is automatically irreducible, two local systems $\mathbb L$ and $\mathbb L'$ are equivalent if and only if there are differed by a character of the absolute Galois group of $\mathbb F_{q^{2n}}$.

Composing the above functors together we obtain an injective functor
$$ M^{gr{1\over 2}}_{Hig\lambda} (\mathbb F_{q^{2n}} ) \simeq M^{FF\,{1\over 2}}_{dR\,\lambda}(\mathbb Z_{q^{2n}})/ \chi \hookrightarrow 
 M^{{1\over 2}}_{\ell-\text{adic}\,\lambda}(\mathbb F_{q^{2n}})/\sim.$$
Over any field $E$ we have 
 $$ M^{gr{1\over 2}}_{Hig\lambda}\otimes E\simeq \mathbb P^1_E,$$ 
 in particular, $$\# M^{gr{1\over 2}}_{Hig\lambda}(\mathbb F_{q^{2n}})=\#\mathbb P^1(\mathbb F_{q^{2n}})=q^{2n}+1.$$
On the other hand by applying a formula due to Hongjie Yu \cite{Y}, which solved a conjecture by Deligne on counting the number of $\ell$-adic local systems on a punctured smooth projective curve $(C,D)/\mathbb F_q$ in terms of the number of semistable parabolic graded Higgs bundles on $(C, D)/\mathbb F_q$, we find $$\#M^{{1\over 2}}_{\ell-\text{adic}\,\lambda}(\mathbb F_{q^{2n}})=q^{2n}+1.$$
This equality implies that the above injective functor is, in fact, bijective
\begin{equation} \label{equ:121} 
M^{gr{1\over 2}}_{Hig\lambda} (\mathbb F_{q^{2n}} ) \simeq M^{FF\,{1\over 2}}_{dR\,\lambda}(\mathbb Z_{q^{2n}})/ \chi \simeq
M^{ {1\over 2}}_{\ell-\text{adic}\,\lambda}(\mathbb F_{q^{2n}})/\sim.	
\end{equation}

\subsection{\bf Constructing abelian scheme over $\mathbb F_q$ via Langlands correspondence and lifting Hodge filtration of relative differential 1-forms to characteristic zero.}

Given a local system $\mathbb L\in M^{ {1\over 2}}_{\ell-adic\,\lambda}(\mathbb F_{q^{2n}}) $ with cyclotomic determinant. Then one shows that the restriction of $\mathbb L$ to the geometric fundamental group is irreducible with infinite local monodromy at least on one puncture. By applying Drinfeld's theorem \cite{D} to $\mathbb L$ there exists an abelian scheme
\[f: A \to \mathbb P^1_{\mathbb F_{q^{2n}}}\]with the bad reduction $\Delta $ over $\{0,1,\infty,\lambda \}$ and such that:
\begin{enumerate}
	\item The abelian scheme is of $\mathrm {GL}_2(K)$-type, where $K$ is the number field generated by traces of Frobenius
	on $(\mathbb L)_x$ at all close points $x\in \mathbb P^1-D$.
	\item Let 
	\[\mathbb V:= R^1_\mathrm{et}f_* \overline{\mathbb Q}_{\ell \, A^0}=
	\bigoplus_{i=1}^g \mathbb L_i\]
	be the $K$-eigen decomposition of the $\ell$-adic local system attached to the family. Then the local system $\mathbb L$ is isomorphic to a $K$-eigen sheaf, say $\mathbb L_1$, and all eigen sheaves have the cyclotomic character. Moreover, the local monodromies matrices of eigen sheaves at any puncture $x \in \{0,1,\infty,\lambda\}$ are of the same type.
	\item Consider the realization of the $\mathcal O_K$-log Dieudonn\'e crystal attached to $f$ over $(\mathbb P^1, \{0,1,\infty,\lambda \})_{\mathbb Z_{q^{2n}}}$, one gets
	\[(V,\nabla,\Phi,\mathcal V)= \Big(R^1_\mathrm{cry}f_* \mathcal O_{A,crys}\Big)(\mathbb P^1_{\mathbb Z_{q^{2n}}}).\]
	The $\mathrm{GL}_2$-structure induces a $K$-eigen sheaves decomposition of $F$-isocrystals
	$$ (V,\nabla,\Phi)\otimes \mathbb Q_p=\bigoplus _{i=1}^g(V,\nabla,\Phi)_{i\,\mathbb Q_p}$$
	where $(V,\nabla,\Phi)_{i,\mathbb Q_p}$ is a $\sigma^f$-log $F$-isocrystal with determinant $\mathcal E_{cy}$. Under $p$-$\ell$ companion, the eigen sheaves correspond to local systems with cyclotomic determinant. The parabolic structure induced by the residue of the connection is of type-$(1/2)$.
\end{enumerate}	
	
 According the bijection 
\eqref{equ:121}, each $(V,\nabla,\Phi)_{i\mathbb Q_p}$ has an integral lattice, which underlies a Fontaine-Faltings module $(V,\nabla, Fil,\Phi)^{FF}_i$ over $\mathbb Z_{q^{2n}}.$ 
Consequently the multiplication field $K$ is unramified at each place above $p$. This shows that the above decomposition can be extended over $\mathbb Z_{q^{2n}}$. In other words, there exists $(V,\nabla,\Phi)_i$ such that 
$$(V,\nabla,\Phi)=\bigoplus_{i=1}^g(V,\nabla,\Phi)_i$$
and 
$$(V,\nabla,\Phi)_i\otimes \mathbb Q_p\simeq 
(V,\nabla,\Phi)^{FF}_i\otimes \mathbb Q_p.$$
Consider the new $\mathcal O_K$-lattice 
$$(V,\nabla,\Phi )^{FF}:=\bigoplus_{i=1}^g(V,\nabla,\Phi)_i^{FF}$$
of $(V,\nabla,\Phi )\otimes \mathbb Q_p$ defined from the Fontaine-Faltings module $$ (V,\nabla, Fil,\Phi)^{FF}:= \bigoplus_{i=1}^g(V,\nabla, Fil,\Phi)_i^{FF}.$$

By extending the coefficient, one gets a Verschiebung on $(V,\nabla,\Phi)\otimes \mathbb Q_p$. By restricting onto the new lattice, one gets a Verschiebung structure $\mathcal V$ on $(V,\nabla, Fil,\Phi)^{FF}$.

According equivalence of the Dieudonn\'e functor, from the new modified $\mathcal O_K$-log Dieudonn\'e crystal $(V,\nabla,\Phi,\mathcal V)^{FF}$, one gets a $p$-isogeny $f': A'\to \mathbb P^1$ of the original abelian scheme such that $(V,\nabla,\Phi,\mathcal V)^{FF}$ is the Dieudonn\'e crystal attached to $f'$.

From the family $f'$, by taking relative differential $1$-forms one gets the natural Hodge filtration on $(V,\nabla)^{FF}\otimes \mathbb F_{q^{2n}}$ given by 
$$E^{1,0}:= R^0f'_*\Omega^1 _{A'/\mathbb P^1} (\log \Delta) \subset (V,\nabla)^{FF}\otimes \mathbb F_{q^{2n}}=R^1_{dR}f'_*(\Omega^\bullet_{A'/\mathbb P^1}(\log \Delta),d),$$
which is a rank-$g$ sub bundle. Since the relative Frobenius $\Phi$ in the Fontaine-Faltings module satisfies the strong $p$-divisible condition with respect to the filtration $Fil$, the Hodge filtration is coincide with the modulo $p$ reduction of the filtration in the Fontaine-Faltings module. In other words, $Fil$ is a filtration lifts the Hodge filtration relative differential $1$-forms attached to $f'$.

\subsection{\bf Lifting abelian scheme from characteristic $p$ to characteristic zero by Grothendieck-Messing-Kato logarithmic deformation theorem} \label{subsec:lifting}
By Zarhin's trick, the fiber product 
\[ f^{'(4,4)}: A^{' 4,4}=(A' \times A^{' \, t})^4 \to \mathbb P^1\]
with the induced $\mathcal O^{(4,4)}_K$-multiplication carries an $\mathcal O^{(4,4)}_K$-principle polarization
$$\iota: A^{' 4,4}\simeq A^{' 4,4\, t}.$$
This $\mathcal O^{4,4}_K$-polarization induces an isomorphism between the $\mathcal O^{4,4}_K$-eigen sheaves of the Dieudonn\'e module and its dual
$$\iota^*: (V,\nabla,\Phi,\mathcal V)^{FF} = \bigoplus_{i=1}^{8g}(V,\nabla,\Phi,\mathcal V)_i^{FF}\simeq \bigoplus_{i=1}^{8g}(V,\nabla,\Phi,\mathcal V)_i^{FF\vee} = (V,\nabla,\Phi,\mathcal V)^{FF\vee} ,$$
which carries the Hodge filtration $ Fil$ and $Fil^\vee$ as liftings of the Hodge bundles $E^{1,0}_{f^{' 4,4}}$ and $E^{1,0}_{f^{' 4,4 t}}$. One checks that
$\text{Gr}_{ E^{1,0}_{f^{' 4,4}} }(V,\nabla)^{FF}\otimes \mathbb F_{q^{2n}} $ and $\text{Gr}_{ E^{1,0}_{f^{' 4,4 t}} }(V,\nabla)^{FF\vee}\otimes \mathbb F_{q^{2n}} $ 
are parabolic stable with respect to the $\mathcal O^{4,4}_K$-eigen sheaves decomposition. Hence, we obtain
\begin{equation} \label{equ:1.4.1}
	\iota^*Fil=Fil^\vee
\end{equation}
By Faltings-Chai Theorem \cite{FC} on the fine arithmetic moduli space $\mathcal A_{8g, 1, 3}=: \mathcal X^0$ of principle polarized abelian varieties with level-3 exists over $\mathbb Z[e^{{2i \pi \over 3}},\, 1/3]$, which is smooth and carries
an universal abelian scheme
\[\mathcal A^0 \to \mathcal X^0.\]
Further more, there exists a smooth Toroidal compactification $\mathcal X \supset \mathcal X^0$ over $\mathbb Z[e^{{2i \pi \over 3}},\, 1/3]$ with a smooth compactification of the universal abelian scheme
\[f: \mathcal A \to \mathcal X\]
such that $\mathcal A \setminus \mathcal A^0$ is a relative normal crossing divisor over $\mathcal X \setminus \mathcal X^0=: \infty$.
For an $\mathcal O^{(4,4)}_E$-principle polarized abelian scheme $f^{' (4,4)}$ together with a lifting of the base $(\mathbb P^1, D)_{\mathbb Z_{q^{2n}}} $\\
 In order to obtain a period map into the fine moduli space $\mathcal X^0$ for $f^{' (4,4)}$ we take the base change
\[\chi: (\mathbb P^{1}, D)_\chi \to (\mathbb P^1, D ) \]
defined by torsion-6 subgroup of $f^{s(4,4)}$. The fiber product of the base change
\[f^{'4,4}_\chi: A^{'4,4}_\chi \to \mathbb P^1_\chi\]
has semistable reduction on $D_\chi$ and carries a level-6 structure. Hence, $f^{s(\lambda4,4)0}_\chi$ induces a log period map
\[\psi: (\mathbb P^{1}, D)_\chi \to (\mathcal X,\infty) \otimes k.\]
By Kato's log deformation theorem \cite{Kat96}, the local liftings of $\psi$ on $W_n(k)$ defines an obstruction-cocycle in $H^1(\mathbb P^1_\chi \otimes k,\psi^* \Theta_{\mathcal X/k}(\log \infty) )$. The local liftings can be glue if and only if the obstruction class vanishes.

By Faltings-Chai \cite{FC}, the universal Kodaira-Spencer map induces an identification
\[\Theta_{\mathcal X}(\log \infty) \simeq S^2E_{\mathcal X}^{0,1}\]
and we may identify the obstruction for lifting of $\psi$ with the obstruction for lifting the Hodge bundle 
$E^{1,0}_{f^{' 4,4}}$ satisfying \eqref{equ:1.4.1}. As $E^{1,0}_{f^{' 4,4}}$ lifts
 we show $\psi$ lifts. Which corresponds to a lifting of the abelian scheme
\[ (f^{'4,4}: A^{' 4,4}_\chi \to \mathbb P^1_\chi)_{\mathbb Z_{q^{2n}}}\]
 with a lifting of $\mathcal O^{4,4}_K$-action, as the the Dieudonn\'e crystal with
 $Fil$ carries an $\mathcal O^{4,4}_K$-action.
 
As $f^{'4,4}$ over $\mathbb P^1_\chi/\mathbb F_{q^{2n}}$ descends to $\mathbb P^1/\mathbb F_{q^{2n}}$ and the Fontaine-Faltings module attached to $f^{' 4,4}$ over $\mathbb P^1_ \chi/W(\mathbb F_{q^{2n}})$ descends to $\mathbb P^1/W(\mathbb F_{q^{2n}})$ we obtain $f^{'4,4}$ descends to $\mathbb P^1/W(\mathbb F_{q^{2n}})$.

\section{\bf Sketch the proof of \autoref{thm:torsion}.}

Let $g: C\to \mathbb P^1$ be the Legendre family. Then one may identify the smooth locus of $g$ with $M_{0,4}$, the moduli space of projective line with $4$-punctures, which sends $\lambda$ to the projective line with punctures at $\{0,1,\lambda,\infty\}$. For any $\lambda\neq 0,1,\infty$, the fiber of $g$ at $\lambda$ is just the elliptic curve given by the double cover $\pi_\lambda: C_\lambda\to \mathbb P^1$ ramified on $\{0,1,\lambda,\infty\}$

For a $\lambda_0 \in M_{0,4}$, take a Higgs bundle $(E,\theta)$ in $M^{gr {1\over 2}}_{Hig\, \lambda_0}.$ Then \autoref{thm:torsion} claims that
$(E,\theta)$ is a motivic Higgs bundle if and only if $(\theta)_0\in \mathbb P^1$ is a torsion point with respect to $\lambda_0$ (i.e. the preimage in 
$\pi_{\lambda_0}^{-1}(\theta)_0\in C_{\lambda_0}$ are torsion points). In the following, we give a sketch of the proof of \autoref{thm:torsion}.

Assume $(E,\theta)$ is motivic. Then the modulo $\mathfrak p$ reduction of $(E,\theta)$ is periodic  for almost all places $\mathfrak p$. According \autoref{thm:LSW}, the modulo $\mathfrak p$ reduction of $(\theta)_0$ is torsion. By a theorem of Pink \cite{Pin04}, it itself is torsion.

Conversely, assume $(\theta)_0$ is a torsion point with order $m$, in the following we show $(E,\theta)$ is motivic. We first show a very special case:

{\bf Case 1.} Assume $\lambda_0$ takes a value in $\mathcal O_K$, the ring of integers of some number field, such that $C_{\lambda_0}$ is an elliptic curve with complex multiplication. Choose a sufficient large place $\mathfrak p$ such that the reduction of $C_{\lambda_0}$ at $\mathfrak p$ is supersingular and $\mathfrak p\nmid m$.

Since the modulo $\mathfrak p$ reduction  $(\bar E, \bar\theta)$ of the Higgs bundle $(E,\theta)$ is also torsion and of order $m$ with $\mathfrak p\nmid m$. By \autoref{thm:LSW}, the reduction $(\bar E, \bar\theta)$ is periodic. According the bijection in \eqref{equ:121}, the modulo $p$ reduction $(\bar E, \bar\theta)$ lifts to a periodic Higgs bundle $(E,\theta)^{per}$ on
$(\mathbb P^1, \{0,1,\lambda_0,\infty\})/W(\mathbb F_{q_0^{2n}})$. In \autoref{subsec:lifting} we show that there exists an abelian scheme $f_{\lambda_0}: A\to \mathbb P^1$ of $\text{GL}_2(K)$-type over characteristic $0$ with bad reduction on
$\{0,\, 1,\, \lambda_0, \, \infty\})$ of type-(1/2) such that $(E,\theta)^{per}$ is an $K$-eigen Higgs bundle attached to $f_{\lambda_0}$. Hence, $(E,\theta)^{per}$ is motivic. In particular $(\theta^{per})_0$ is also torsion and has the some modulo $\mathfrak p$ reduction as $\theta_0$. 

We claim that $(E,\theta)\simeq (E, \theta)^{per}$. In particular, $(E,\theta)$ is motivic. Since $C_{\lambda_0}$ has complex multiplication, the field generated by 
$p$-torsion point must be ramified above $p$. In particular, the order of $(\theta^{per})_0$ does not divided by $p$. Together with the fact that $(\theta)_0\equiv (\theta^{per})_0 \pmod{\mathfrak p}$, one gets  
$(E,\theta)\simeq (E,\theta)^{per}$.

{\bf Case 2.} Let $\Sigma_m\subset C$ be the $m$-torsion (multiple) section, $T_m=\pi(\Sigma_m)\subset M_{0,4}.$
Then $T_m$ is \'etale over $M_{0,4}$.
Let $T'_m$ be the irreducible component of $T_m$ containing $(\theta)_0$.

Recall the fact that the set of elliptic curves with complex multiplication are dense in the moduli space. We can find a subset of infinitely many $\{\lambda_i\}$ with the some modulo $\mathfrak{p}$-reduction such that all $C_{\lambda_i}$ are elliptic curve with complex multiplication. Choose a point $z_i$ in the intersection of $T'_m$ and the fiber above $\lambda_i$ for each $i$.

With the method in case 1, we find abelian schemes $f_{\lambda_i}: A_{\lambda_i}\to \mathbb P^1$ of $\text{GL}_2(K)$-type with bad reduction on $\{0,\,1,\,\lambda_i,\, \infty \},$ such that $(E,\theta)_{z_i}$ is an $K$-eigen Higgs bundle attached to $f_{\lambda_i}$. One shows that infinitely many of $\{f_{\lambda_i}\}$ glue together into an abelian scheme $f: A\to M_{0,4}$. Consequently, the subset $\{z_i\}$, where $z_i=(\theta_{z_i})_0 $, glue into a component $Z_m$ of $\Sigma_m$ such that for any $z_\lambda \in Z_m\cap (\mathbb P^1, \{ 0,\,1,\,\lambda,\, \infty  \}) $ the Higgs bundle $(E,\theta)_{z_\lambda}$ with $(\theta_{z_\lambda})_0=z_\lambda$ is a motivic Higgs bundle on $(\mathbb P^1, \{ 0,\,1,\,\lambda,\, \infty\}).$ 

By the construction of $Z_m$ and $T'_m$, Both of then are irreducible, relative dimensional one over the base and their intersect set is infinite. Hence $Z_m =T_m$. In particular $(E,\theta)$ is motivic.

\begin{acknowledgement*}
The authors warmly thank Raju Krishnamoorthy, Mao Sheng, \\ Carlos Simpson, Rui-Ran Sun, Hong-Jie Yu and Shing-Tung Yau for helpful discussions.
\end{acknowledgement*}

\end{document}